\documentclass[12pt]{amsart}
\usepackage{amsmath,amssymb, euscript, graphicx, hyperref}

\usepackage{latexsym,array,delarray,amsthm,amssymb,epsfig,setspace,tikz,amsmath,enumerate,mathrsfs,graphicx, mathtools,cite}
\usepackage{geometry}
\usepackage{tabularx}

\newtheorem{thm}{Theorem}[section]
\newtheorem{lem}[thm]{Lemma}
\newtheorem{prop}[thm]{Proposition}
\theoremstyle{definition}
\newtheorem{defn}[thm]{Definition}

\newcommand{\blackboard}[1]{\ensuremath{\mathbb{#1}}}

\newcommand{\Z}{\blackboard{Z}}

\newcommand{\R}{\blackboard{R}}

\newcommand{\F}{\ensuremath{\mathbf{F}}}

\begin{document}

\address{Azer Akhmedov, Department of Mathematics,
North Dakota State University,
Fargo, ND, 58102, USA}
\email{azer.akhmedov@ndsu.edu}

\address{James Thorne, Department of Mathematics, Computer Science, and Statistics,
Gustavus Adolphus College,
800 West College Avenue,
Saint Peter, MN 56082 
 AND
North Dakota State University,
Fargo, ND, 58102, USA}
\email{james.thorne@ndsu.edu}

 \begin{center} {\bf \Large Non-bi-orderable one-relator groups without generalized torsion} \end{center}
 
 \bigskip
 
 \begin{center} {\bf Azer Akhmedov and James Thorne} \end{center}
 
 \medskip
 
 \begin{center} {\small ABSTRACT: We construct examples of non-bi-orderable one-relator groups without generalized torsion. This answers a question asked in \cite{CGW}.}
 \end{center}

  \bigskip

\section{Introduction}

 A left-order on a group is a linear order which is invariant under the left translations, and a bi-order is a linear order which is invariant under both left and right translations. Bi-orderable groups turn out to be, perhaps not too narrow, but still a rather special sub-class in the class of all left-orderable groups.   

 \medskip 
 
  The existence of a torsion element is clearly an obstruction to left-orderability of a group. Similarly, existence of a generalized torsion element is an obstruction to bi-orderability. There exist examples of groups without a torsion which are not left-orderable, and examples of non-bi-orderable groups without a generalized torsion also exist (\cite{MR}, \cite{CR}).  On the other hand, for many interesting classes of groups, the absence of a torsion element turns out to be equivalent to left-orderability. One of the non-trivial examples is the class of one-relator groups (\cite{B}, \cite{H1}, \cite{H2}, \cite{CR}). Interestingly, torsion-freeness of a one-relator group can be immediately detected from the relator, namely, a one-relator group on generators $a, b, \dots , c$ contains a torsion element iff the relation of a group is in the form $r^n$  for some $n\ge 2$ and $r = r(a,b,\dots , c)\in F$ where $F$ is the free group on generators $a,b, \dots , c$ (see \cite{LS}). 

 \medskip 

 In \cite{CGW}, the authors study bi-orderability of one-relator groups and prove that 1) if a one-relator group with the so-called tidy relation is 
 bi-orderable, then its Alexander polynomial has at least one positive real root; 2) if all the roots of Alexander polynomial of a one-relator group with monic relation are real and positive, then the group is bi-orderble. They also ask (See Question 3) if every one-relator group without a generalized torsion is bi-orderable.  

 \medskip

  In this note, we produce non-bi-orderable examples of one-relator groups which do not possess a generalized torsion. Let us recall that for a group $G$, an element $g\in G\backslash \{1\}$ is called {\em a generalized torsion} if there exist $x_1, \dots , x_n\in G$ for some $n\geq 1$  such that $(x_1gx_1^{-1})\dots (x_ngx_n^{-1}) = 1$. 

\medskip 

  The following theorem produces an explicit example of a one-relator non-bi-orderable group without a generalized torsion. 

   \medskip 
   
 \begin{thm} \label{thm:special} For all $q\geq 8$, the one-relator group $$\Gamma = \langle t, a \ | \ a^{2}ta^{-1}t^{-1}[a, tat^{-1}]^qt^2a^2t^{-2}  \rangle$$ is not bi-orderable and does not contain a generalized torsion. 
   \end{thm}

 \bigskip 

  Let us notice that the Alexander polynomial (\cite{CGW}) of $\Gamma $ is $2t^2-t+2$ (so none of the roots are real and positive). Notice that our group is a 2-generated 1-relator group, thus we also answer the specific question asked in \cite{CC}. The group that we present can be generalized in various directions. We do not make an attempt here for this, and in particular, we make no effort to clarify how small the exponent $q$ can be taken. 

  \bigskip 

  In the class of 3-manifold groups, it is conjectured in \cite{MT} that the existence of a generalized torsion is equivalent to non-bi-orderability. The authors of \cite{MT} verify it for all non-hyperbolic, geometric 3-manifolds (in particular, for Seifert fibered manifolds and for Sol manifolds) and in \cite{IMT}, further evidence for the conjecture are provided. Recently, the conjecture has settled negatively in \cite{CC}. In the same work, a computation is provided to show that the example in the previous version of our paper has a generalized torsion. In the current version, we replace the flawed argument from the previous version with a very short and efficient argument.    

  \medskip 

   For any class $\mathcal{C}$ of groups, we say $\mathcal{C}$ satisfies the property $T_O$ if for any group $G$ in this class, $G$ is left-orderable iff $G$ is torsion-free. Similarly, we say $\mathcal{C}$ satisfies property $T_{BO}$ if for any group $G$ in $\mathcal{C}$, $G$ is bi-orderable iff $G$ has no generalized torsion. For the class of 3-manifold groups, both properties $T_O$ and $T_{BO}$ fail (see \cite{G}, \cite{CC}); and for the class of one-relator groups, property $T_O$ holds, whereas property $T_{BO}$ fails. We summarize this in the following table:
   
   \bigskip 
   
    \begin{tabularx}{1.0\textwidth} { 
  | >{\raggedright\arraybackslash}X 
  | >{\centering\arraybackslash}X 
  | >{\raggedleft\arraybackslash}X | }
\hline
\textbf{Groups} & \textbf{$T_O$} & \textbf{$T_{BO}$} \\
\hline
One-relator Groups & holds \cite{B}, \cite{H2} & fails (current paper and \cite{CC})\\
\hline

3-Manifold Groups & fails \cite{G}  & fails in general \cite{CC}, holds in some cases \cite{IMT} \\
\hline

\end{tabularx}

    \vspace{1cm} 

    \section{Proof of Theorem \ref{thm:special} }

\medskip

 \medskip 

    A one-relator group $G$ in $n\geq 3$ generators can always be embedded into a group $H$ with two generators and one relator. For the latter, by a sequence of Shreier moves we can obtain a one-relator presentation $\langle t, h \ | r\rangle $ of it such that sum of exponents of one of the generators (say, of the  generator $t$) is zero. This allows us to view $H$ as a semidirect product $\Z \ltimes H_0$ where the acting group $\Z$ is generated by $t$, the normal subgroup $H_0$ is generated by elements $h_n = t^nht^{-n}, n\in \Z$ and the $\Z$ action is given by the conjugation $th_nt^{-1} = h_{n+1}, n\in \Z$. 

    \medskip 

    The relation $r$ can be written as a reduced word $R(h_0, \dots , h_k)$ in generators $h_0, \dots , h_k$ for some $k\geq 0$ and this allows to view $H$ as an HNN extension $(H, A, B, t)$ where $$A = \langle h_0, \dots , h_k \ | \ R(h_0, \dots , h_k) \rangle , \ \mathrm{and} \ B = \langle h_1, \dots , h_{k+1} \ | \ R(h_1, \dots , h_{k+1}) \rangle $$
 where one can even claim that the word $R$ is shorter than the word $r$ in length. This inductive phenomenon together with the HNN extension structure become very helpful in understanding one-relator groups. One of the strong tools in utilising the HNN extension structure is given by Britton's Lemma. In our approach, we will be using the amalgamated free product structure of one-relator groups (rather than HNN extension) since the group $H_0$ above is a direct limit of groups $H_{m,n} = \langle h_m, \dots , h_n \rangle , m\leq n$ and for all $m < n$ with $n-m \geq k$, we can write $H_{m,n}$ as the amalgamated free product $H_{m,n-1}\ast _{H_{m+1,n-1}}H_{m+1,n}$. (In our proof, we will use a slightly different but related amalgamated free product decomposition of the groups $H_{m,n}$.) 

  \medskip 
  
  The following result is well known (see \cite{LS}); it can be viewed as an analog of Britton's Lemma and plays an important role in understanding amalgamated free products similar to the role that Britton's Lemma plays in understanding HNN extensions.

  \medskip 

  \begin{prop} \label{thm:amalgam} Let $A, B, C$ be groups with monomorphisms $\phi :C\to A$ and $\psi :C\to B$ and $G = A*_CB$ be the amalgamated free product with respect to these monomorphisms. Let also $T_A = A\backslash C, T_B = B\backslash C$ and $w = g_0g_1\dots g_n, n \geq 1$ such that for all $1\leq  i \leq n$  
 
 (i) if $g_{i-1}\in A$, then $g_i\in T_B$
 
 (ii) if  $g_{i-1}\in B$, then $g_i\in T_A$
 
 Then $w\neq 1\in  A*_CB$.

  \end{prop}

  \medskip 

   We will need another technical statement,  which we will formulate in the following as a lemma. First, we need to introduce some useful notation. 

   \begin{defn} Let $G$ be a group, $F\trianglelefteq G$ be a free group on a countable subset $S$, $\pi :F\to H$ be an epimorphism, and let $|z|$ denotes the length of $z\in F$ in the left-invariant Cayley metric with respect to the generating set $S$ of $F$.

   \medskip 
   
   i) For all $x\in F\backslash \ker (\pi )$, we let $C(x) = \{g\in G : g \ \mathrm{ is \ a \ product \ of \ conjugates \ of} \ x \ \mathrm{in} \ G\}$, and $M(x) = \{\zeta \Theta \eta^{-1} : \Theta \in C(x), \zeta , \eta  \in F, \pi (\zeta ) = \pi(\eta )\}$; 

   \medskip 
   
   (ii) For all $x\in F\backslash \ker (\pi ), \Theta \in C(x)$, we let $M(x, \Theta ) = \{\zeta \Theta \eta^{-1} : \zeta, \eta  \in F, \pi (\zeta ) = \pi(\eta )\}$; 
   
   \medskip

   (iii) For all $x\in F\backslash \ker (\pi ), W \in M(x)$, fixing $x$, we also let $$L(W) = \displaystyle \mathop{\min }_{\displaystyle \mathop{\Theta \in C(x), W\in M(x,\Theta )}} |\Theta | .$$

   \end{defn}

\medskip 

    Now we are ready to state the lemma. Its essence is indeed simple, but in its technical statement, we have tried to capture all the details of the situation in which it applies. 

    \medskip

   \begin{lem} \label{prop:shuffling} Let $G$ be a group, $F\trianglelefteq G$ be a free group on a countable subset $S$, $\pi _1 :F\to K, \pi _2:K\to H$ be epimorphisms, $\pi = \pi _2\circ \pi _1$,  and let $|z|, |z'|_K$ denote the lengths of $z\in F$ and $z'\in K$ in the left-invariant Cayley metrics with respect to the generating sets $S$ of $F$ and $\pi _1(S)$ of $K$ respectively. Let $x\in F\backslash \ker (\pi )$ and  $C(x)\cap ker \pi \neq \emptyset $ such that $|x|$ is minimal. Let also $W_0, W_1, W_2, W, \Theta \in F, \Theta  = W_1W_0W_2$ with $\pi (W_0) = \pi (W) = 1\in H$ and $|\Theta | = |W_1|+|W_0|+|W_2|$ such that for any subword $V$ of the reduced word $\Theta $ in $F$, we have $|V|_K = |V|$, moreover, $\Theta \in C(x), W\in M(x, \Theta )$ such that $|\Theta |$ is minimal (so $L(W) = |\Theta |$). Then there exist $A, B, C, D, \Theta ', W'\in F$ such that $\Theta '\in C(x), W'\in M(x, \Theta ')$, $|\Theta '|\leq |\Theta |$ (so $L(W')\leq L(W)$), and $\Theta ' = ABCD$ where either $A$ is a suffix of $W_1$, $B$ is a prefix of $W_2$, $C$ is a subword of $W_0$, $D$ is a suffix of $W_2$ disjoint from $B$, or $A$ is a prefix of $W_1$, $B$ is a subword of $W_0$, $C$ is a suffix of $W_1$ disjoint from $A$, $D$ is a prefix of $W_2$.

    \end{lem}

    \begin{proof}  By the minimality assumption of $|x|$, and since $\pi (W_0) = 1$, $W_0$ cannot be a subword of $x$. 

\medskip 

  For reduced words $V_1, \dots , V_m$ in the free group $F$ on the set $S$, let $\gamma $ be the path (geodesic) that represents $V := V_1\dots V_m$ as a reduced word in the Cayley graph of $F$ with respect to the generating set $S$. We will say that a subword $v$ of $V_1\dots V_m$ {\em lies in the image of} $V_i$ {\em in } $V$ for some $1\leq i\leq m$, if the geodesic corresponding to $v$ is contained in the intersection of the geodesic $\gamma $ and the geodesic $\beta _i$ from $V_1\dots V_{i-1}$ to $V_1\dots V_{i}$. (Thus, $v$ is a subword of $V_i$.)
  
\medskip

     Let $\Theta = (g_1xg_1^{-1})\dots (g_nxg_n^{-1})$. We may assume that $|g_ix| = |g_i|+|x| = |xg_i^{-1}|$ for all $1\leq i\leq n$ (otherwise, we can proceed by replacing $x$ in each factor with a cyclic permutation of it). Let $W_0$ be a subword of $g_ixg_i^{-1}$ for some $1\leq i\leq n$ such that $W_0$ lies in the image of $g_ixg_i^{-1}$ in $W$. If $W_0 = f_1xf_2^{-1}$ for some subwords $f_1, f_2$ of $g_i$, then let $g_i = h_1f_1 = h_2f_2$ with $|g_i| = |h_1|+|f_1| = |h_2|+|f_2|$. Let $|f_1| \geq |f_2|$ (the case $|f_1| \leq |f_2|$ is similar). Then $h_1$ is a prefix of $h_2$, so let $h_2 = h_1h$. Then we replace $g_i$ with $h_2x$; then, since $|f_1xf_2^{-1}| > |hx|$, we obtain a contradiction to the minimality of $L(W)$. 
     
     \medskip 
     
     Then, without loss of generality, we may assume that $g_i = f_1f_2, x = x_1x_2$ and $W_0 = f_2x_1$ with $|W_0| = |f_2|+|x_1|$. If $|f_2| < |x_1|$, then by replacing $x = x_1x_2$ with $f_2^{-1}x_2$ we obtain a contradiction to the minimality of $|x|$. So we may assume that $|f_2| \geq |x_1|$. Then, replacing $f_2$ with $x_1^{-1}$ we obtain a new word $\Theta ' = U_1f_1x_2x_1f_1^{-1}U_2$ where $W_1 = U_1f_1, W_2 = x_2f_2^{-1}f_1^{-1}U_2$ and $W_0 = f_2x_1$ (so $\Theta = W_1W_0W_2 = U_1f_1f_2x_1x_2f_2^{-1}f_1^{-1}U_2 $). Then we let $A = W_1, B = x_2, C = x_1, D = f_1^{-1}U_2$.

    \medskip

    Now, assume that $W_0$ does not lie in the image of any $g_ixg_i^{-1}, 1\leq i\leq m$. Then, without loss of generality, we may assume that $\Theta  = W_1ugg^{-1}vW_2$ where $W_0 = uv$, and $W_1ug$ and $g^{-1}vW_2$ are products of conjugates of $x$ in $G$. Then $ugW_1$ and $W_2g^{-1}v $ are products of conjugates of $x$ in $G$, and so is $ugW_1W_2g^{-1}v$.  On the other hand, $\pi (u) = \pi (v)^{-1}$ therefore $\pi (ug) = (\pi (g^{-1}v))^{-1}$. But $|W_1W_2| < |\Theta |$ which contradicts the minimality of $L(W)$.

    \end{proof}
    
    \bigskip 

  We are now ready to start the proof of Theorem \ref{thm:special}. First, we prove that $\Gamma $ is not bi-orderable.     

\medskip 

    Let $a_n = t^nat^{-n}, n\in \Z$. Then we have $a_0^2a_1^{-1}[a_0,a_1]^{q}a_2^2 = 1$. 

    \medskip

     Now, assume that $<$ is a bi-order on $\Gamma $. Without loss of generality, we may assume that $a_0 > 1$. Then $a_n > 1$ for all $n\in \Z$. If $a_0$ and $a_1$ are comparable, then $[a_0,a_1]$ is infinitesimal with respect to both, hence $a_1^{-1}[a_0,a_1]^{q} > a_1^{-2}$, hence $a_0^2a_1^{-1}[a_0,a_1]^{q}a_2^2 > a_0^2a_1^{-2}a_2^2 > 1$; contradiction. Thus, we may and will assume that $a_0$ and $a_1$ are not comparable.

     \medskip 
     
     If $a_1 > a_0$, since the positive cone is invariant under a conjugation, we obtain $1 < a_0 < a_1 < a_2$; and if $a_1 < a_0$, then $1 < a_2 < a_1 < a_0$.  In the former case, $a_0$ is infinitesimal with respect to $a_1$. Then $a_1$ (hence also $a_0$) is infinitesimal with respect to $a_2$; therefore $a_0^2a_1^{-1}[a_0,a_1]^{q}$ is infinitesimal with respect to $a_2$, hence $a_0^2a_1^{-1}[a_0,a_1]^{q}a_2 > 1$, hence $a_0^2a_1^{-1}[a_0,a_1]^{q}a_2^2 > 1$ which is again a contradiction.
     In the latter case, $a_1$ is infinitesimal with respect to $a_0$; then, $[a_0, a_1]$ and $a_2$ are also infinitesimal with respect to $a_0$, and we obtain $a_0^2a_1^{-1}[a_0,a_1]^{q}a_2^2 > a_0a_1^{-1}[a_0,a_1]^{q}a_2^2 > 1$; contradiction.

 \medskip 

 Now, assume that $\Gamma $ contains a generalized torsion $\tau $. 

 \medskip 
 
  Let $H$ be the subgroup of $\Gamma $ generated by all $a_n, n\in \Z$. Then $H$ is a normal subgroup and $\Gamma \cong \Z\ltimes H$ where $\Z$ acts on $H$ by shift, i.e. by the $t$-conjugation $ta_nt^{-1} = a_{n+1}, n\in \Z$. Then necessarily $\tau \in H$.

  \medskip 

  Elements of $H$ can be written (not uniquely) as reduced words in the alphabet $\mathcal{A} = \{a_n , n\in \Z\}$. Then $H$ can be given by the presentation $$H = \langle a_n, n\in \Z \ | \ a_{n}^2a_{n+1}^{-1}[a_n, a_{n+1}]^qa_{n+2}^2 = 1,  n\in \Z \rangle .$$ 

  \medskip 

   Now, for all $n\geq 0$, let $G_n$ be the subgroup of $H$ generated by $a_0, \dots , a_n$ and for all $m \leq n$, let $H_{m,n}$ be the subgroup of $H$ generated by $a_m, \dots , a_n$. Then $H_{m,n}\cong G_{n-m}$ by the obvious isomorphism sending $a_{m+j}$ to $a_j, 0\leq j\leq n-m$.  Furthermore, $$G_0 \cong \Z, G_1\cong \mathbb{F}_2, G_2 \cong \langle a,b,c \ | \ a^{2}b^{-1}[a,b]^qc^2\rangle .$$ Additionally, for all $n\geq 3$ we make the following key observations about $G_n$: The group $G_n$ is isomorphic to an amalgamated free product $G_{n-1}\ast _{\Z} \Z$, more precisely, we have $$G_n \cong G_{n-1}\ast _{\Z }H_{n,n} \cong G_{n-1}\ast _{\Z} \Z$$  where the amalgamation is over the cyclic group $\Z = \langle t \rangle $, the monomorphism $\phi _1:\Z\to G_{n-1}$ is just the map given by $\phi _1(t) = a_{n-2}^2a_{n-1}^{-1}[a_{n-2}, a_{n-1}]^{q}$, and the monomorphism $\phi _2:\Z\to H_{n,n}$ is given by $\phi _2(t) = a_{n}^{-2}$. In addition, we also let $$K_n = \mathcal{F}/\langle a_{n-2}^2a_{n-1}^{-1}[a_{n-2}, a_{n-1}]^{q}a_{n}^2 \rangle , n\in \Z $$ 
    
where $\mathcal{F}$ is the free group formally generated by the alphabet $\mathcal{A} $.  $\Z$ acts on the alphabet $\mathcal{A} $ by the shift $a_n\to a_{n+1}, n\in \Z$ hence also on the groups $\mathcal{F}$ and $H$ as automorphisms. We let $\pi_1:\mathcal{F}\to K_n, \pi :\mathcal{F}\to H $ be the quotient maps and $\pi _2:K_n\to H$ be the epimorphism such that $\pi _2\circ \pi_1 = \pi $. Abusing the notation, we will write $twt^{-1}$ for the shift of a word $w\in \mathcal{F}$. In addition, we will write $|w|, |w|_{K_n}$ to denote the lengths of $w\in \mathcal{F}$ in the Cayley graphs of $\mathcal{F}$ and $K_n$ respectively, with respect to the generating set $\mathcal{A}$.

   \medskip 

   For a reduced non-empty word $w\in \mathcal{F}$ in the alphabet $\mathcal{A}$, we call the generator $a_n$ {\em dominant} if $n = \max \{i\in \Z : a_i \ \mathrm{occurs \ in} \ w\}.$ We also say that $a_n$ occurs at $k$ places in $w$, if as a reduced word, $w$ equals $w_0a_n^{j_1}w_1\dots w_{k-1}a_n^{j_k}w_k$ where $w_i, 0\leq i\leq k$ are words in the alphabet $\{a_i : i < n\}$ (i.e. $w_i$ do not contain $a_n^{\pm 1}$), $w_i, 1\leq i\leq k-1$ are non-empty words and $j_1, \dots , j_k$ are non-zero integers.   
   
   \medskip 
   
   A reduced non-empty word $w$ will be called {\em suitable} if its dominant generator has occurrence with an odd exponent (i.e. not all exponents of it are even), moreover, $w$, as a reduced word,  is not in form $w = v_0xv_1yv_2$ where $v_1$ is a non-empty subword representing the identity element in $H$ and $x,y\in \{a_n, a_n^{-1}\}$ with $a_n$ being the dominant generator.  $w$ will be called {\em strongly suitable} if, in addition to being suitable,  all the exponents of its dominant generator are odd. If $w = w_0a_n^{j_1}w_1\dots w_{k-1}a_n^{j_k}w_k$ (as described in the previous paragraph) is strongly suitable or even suitable and $k$ is minimal possible, then, in the expression $w = v_0xv_1yv_2$, the subword $v_1$ cannot represent a power of the dominant generator either, therefore, by Proposition \ref{thm:amalgam}, $w$ cannot represent the identity element of $H$. $w$ is called {\em weakly suitable} if as a reduced word it is in the form $w_0a_n^{j_1}w_1\dots w_{s-1}a_n^{j_s}w_s$ where $j_1, \dots , j_s$ are odd integers, $w_i, 1\leq i\leq s-1$ do not represent identity element of $H$,  $w_i, 0\leq i\leq s-1$ do not end with $a_n^{\pm 1}$, $w_i, 1\leq i\leq s$ do not start with $a_n^{\pm 1}$.\footnote{Notice that in the definition of a weakly suitable word, the subwords $w_i, 0\leq i\leq s$ are still allowed to contain $a_n^{\pm 1}$. We just impose conditions on them about starting or ending with  $a_n^{\pm 1}$.}  In addition to our observations about strongly suitable and suitable words, we make an important observation that, by Proposition \ref{thm:amalgam},  if $s$ is minimally possible, then the weakly suitable word $w$ cannot represent the identity element of $H$ either.

   \medskip

   We may assume that $\tau $ is a generalized torsion with the shortest possible length in the group $\mathcal{F}$ in the alphabet $\mathcal{A}$; then it also has the shortest possible length in $H$ in the alphabet $\mathcal{A}$. We have $$(t^{n_1}g_1\tau g_1^{-1}t^{-n_1})(t^{n_2}g_2\tau g_2^{-1}t^{-n_2})\dots (t^{n_r}g_r\tau g_r^{-1}t^{-n_r}) = 1$$ for some $r\geq 1$ where $g_j\in H, n_j\in \Z, 1\leq r$ are chosen such that if $t^{n_i}g_it^{-n_i}, t^{n_i}\tau t^{-n_i}$ are written as reduced words $V_i, \tau _i$ respectively in the alphabet $\mathcal{A} $, then the words $$\Theta \in C(\tau ) \ \mathrm{and} \ W  = V_1\tau _1V_1^{-1}\dots V_r\tau _rV_r^{-1}\in C(\tau )\cap M(\tau , \Theta )$$ viewed as reduced words in the alphabet $\mathcal{A} $ are such that $\Theta $ has minimal length in the group $K_n$ (we will assume that this minimal length is strictly positive) with respect to the alphabet $\{a_i : i\leq n\}$ where $a_n$ is the dominant generator of $\Theta $ and the following conditions hold:

 \medskip 
 
 i) for all $1\leq i\leq r$, the reduced word $\tau _i$ is just a cyclic permutation of letters of $t^{n_i}\tau t^{-n_i}$;

\medskip 

 ii) for all $1\leq i\leq r$, $|V_i\tau _i| = |V_i| + |\tau _i| = |\tau _i(V_i)^{-1}|$. \footnote{Let us emphasize that $$|W| \leq \displaystyle \mathop{\sum } _{i=1}^r|V_i|+|\tau _i|+|V_i^{-1}| = \displaystyle \mathop{\sum }_{i=1}^r(2|V_i|+|\tau |).$$}

   \medskip 
   
   Assume that the dominant generator $a_n$ occurs $k$ times in the reduced word of $\Theta $ where $k$ is minimally possible.\footnote{Notice we first chose the length of $\tau $ to be minimal; then, for those choices of $\tau $, we chose $|\Theta |_{K_n}$ minimal; then, for the choices of $\tau $ and $|\Theta |_{K_n}$, we choose $k$ to be minimal.} We will assume that $k$ is strictly positive. We can write $\Theta $ as $U = W_0a_n^{j_1}W_1a_n^{j_2}\dots W_{k-1}a_n^{j_k}W_k$ in a reduced form in $\mathcal{F}$ where $j_1, \dots , j_k$ are non-zero integers and the letter $a_n$ does not occur in $W_i, 0\leq i\leq k$ such that $|U|_{K_n} = \displaystyle \mathop{\sum }_{0\leq i\leq k}|W_i|_{K_n} + \displaystyle \mathop{\sum }_{1\leq i\leq k}|j_i|$. If $W_i=1$ in $H$ for some $1\leq i\leq k-1$, then by Lemma \ref{prop:shuffling} we can replace $U$ with $U'= W_0'a_n^{i_1}W_1'a_n^{i_2}\dots W_{s-1}'a_n^{i_s}W_s'$ as a reduced word in $\mathcal{F}$ such that $U'$ is also a product of conjugates of $\tau , W\in M(\tau, U')$, $W_i',0\leq i\leq s$ do not contain $a_i$ for all $i\geq n$, $s < k$, and if $a_m$ is the dominant generator of $U'$ with $m < n$, then $|U'|_{K_m} < |U|_{K_n}$. This contradicts the minimality of $|U|_{K_n}$ or minimality of $k$.  Thus, we may assume that $W_i\neq 1$ in $H$, for all $1\leq i\leq l\leq k-1$. By the same argument, we may also assume that for all $1\leq i \leq l\leq k-1$, in the word $U$, $W_ia^{j_{i+1}}W_{i+1}\dots a^{j_l}W_l\neq 1$ in $H$.
   
   \medskip

   Notice that $|a_n^j|_{K_n} = j$, for all $j\in \Z$. Therefore, if $U$ is suitable, then in addition to having $W_i\neq 1, 1\leq i\leq k-1$, (even though the choice of $k$ is after the choice of minimal length of $|U|_{K_n}$) we can also claim that for all $1\leq i\leq k-1$, $W_i$ is not a power of $a_n$. Then $U\neq 1$ in $H$ by Proposition \ref{thm:amalgam}, and we are done. Thus, let us assume that $U$ is not suitable. Then all the exponents of the dominant generator $a_n$ are even and we can replace all the occurrences of $a_n^2$ with $[a_{n-1},a_{n-2}]^qa_{n-1}a_{n-2}^{-2}$. \footnote{That is, if $U = U_1a_n^{2l}U_2$ as a reduced word where $U_1$ does not end with $a_n^{\pm 1}$ and $U_2$ does not start with $a_n^{\pm 1}$, then we replace $U$ with $U_1([a_{n-1},a_{n-2}]^qa_{n-1}a_{n-2}^{-2})^lU_2$}. Notice that $$|[a_{n-2},a_{n-1}]^{q-3}| = 4(q-3) > \frac{1}{2}(4q+5) =  \frac{1}{2} |a_{n-2}^{2}a_{n-1}^{-1}[a_{n-2},a_{n-1}]^qa_{n}^2| .$$ 
   
   Then, in the word $U = W_0a_n^{j_1}W_1a_n^{j_2}\dots W_{k-1}a_n^{j_k}W_k$, if $j_l = -j_{l+1}$ for some $1\leq l\leq k-1$, then by replacing the squares of $a_n$ in $a_n^{j_l}W_la_n^{j_{l+1}}$, by the minimality of $|U|_{K_n}$,  we obtain a reduced word $u_1\omega \Omega \omega ^{-1}u_2$ where $\omega \in \{\eta , \eta ^{-1}\}$ with $\eta = [a_{n-1}, a_{n-2}]^2a_{n-1}$ and $\Omega $, being a reduced word in the alphabet $\{a_i : i < n\}$, does not represent the identity element in $H$. On the other hand, if $j_l = j_{l+1}$ for some $1\leq l\leq k-1$, then by replacing the squares of $a_n$ in $a_n^{j_l}W_la_n^{j_{l+1}}$ we obtain $u_1\omega _1\Omega \omega _2u_2$ with $\omega _1, \omega _2\in \{\eta , \eta ^{-1}, \theta , \theta ^{-1}\}$, where if $\Omega =1$ in $H$, then (since $a_{n-2}^{\delta _1}a_{n-1}^ja_{n-2}^{\delta _2}\neq 1$ in $H$ for all $\delta _1, \delta _2\in \{-1,1\}, j\in \Z\backslash \{0\}$) we obtain a reduced word $u_1'\omega '\Omega '\omega ''u_2'$ where $\omega ', \omega '' \in \{\eta , \eta ^{-1}, \theta, \theta ^{-1}\}$ with $\theta = [a_{n-1}, a_{n-2}]a_{n-1}$ and $\Omega '$ a reduced word in the alphabet $\{a_i : i < n\}$ not representing the identity element in $H$. Thus, again by the minimality of $|U|_{K_n}$, from the word $U$ we obtain a reduced word $V = \Omega _0\omega _1\Omega _1\omega _2\dots \Omega _{k-1}\omega _k\Omega _k$ where $\omega _i\in \{\eta , \eta ^{-1}\}, 1\leq i\leq k$ and $\Omega _i, 0\leq i\leq k$ are reduced words in the alphabet $\{a_i : i < n\}$ such that if $\omega _l = \omega _{l+1}^{-1}$ for some $1\leq l\leq k-1$, then $\Omega _l$ does not represent the identity element in $H$ (if $\omega _l = \omega _{l+1}$, we allow $\Omega _l = 1\in H$). Such a word $V$ is necessarily weakly suitable because we can also write it in the form $V = \Omega _0'\omega _1\Omega _1'\omega _2\dots \Omega _{k-1}'\omega _k\Omega _k'$ where $\omega _i\in \{\eta , \eta ^{-1}, \theta , \theta ^{-1}\}, 1\leq i\leq k$ and $\Omega _l'$ are reduced words in the alphabet $\{a_i : i < n\}$ not representing the identity element in $H$. Then by Proposition \ref{thm:amalgam}, $V\neq 1$ in $\Gamma $. Contradiction.   

   \medskip 

   Now, it remains to show that the minimal length of $\Theta $ viewed as a reduced word in the alphabet $\mathcal{A} $ in the group $K_n$ can be assumed to be positive, moreover, the dominant generator $a_n$ occurs in it. Notice that $\pi (W) = \pi (\Theta ) = 1\in H$ and $\Theta , W\in C(\tau )$, so it suffices to show that the minimal length $|W|_{K_n}$ of $W = V_1\tau _1V_1^{-1}\dots V_r\tau _rV_r^{-1}$  can be assumed to be positive. We will indeed show the positivity of the length and the occurrence of the dominant generator for any word in $C(\tau )\cap ker \pi $. For this we will work with $\tau $ where all occurrences of $a_n^2$ are replaced with $[a_{n-1}, a_{n-2}]^{q}a_{n-1}a_{n-2}^{-2}$ (so $\tau $ does not necessarily have a minimal length). 
   
   \medskip 
   
   Let $\overline{V_i} = V_i\tau _iV_i^{-1}, 1\leq i\leq r$ and we may assume that for all $1\leq i\leq r$, there exists a subword of the reduced word $W$ which lies in the image of $\overline{V_i}$. In addition, we will replace each $\overline{V_i}, 1\leq i\leq r$ with a minimal length representation under the condition that no $a_n^2$ occurs in any of them. To show that the length $|W|_{K_n}$ is positive, it suffices to prove that the image of $W$ is a non-identity element in $K_n$.

   \medskip 
   
   Let the dominant generator $a_n$ occur in $k$ places in $W$ and let $k$ be minimal. If $a_n$ does not occur in $\tau $, then, since $a_n^2$ occurrences are replaced in each $\overline{V_i}, 1\leq i\leq r$, we will find that $a_n$ occurs in $W$ unless $W$ is a product of conjugates in the subgroup $\mathcal{F}_n$ of $\mathcal{F}$ generated by $\{a_i : i < n\}$, but in that case $W\neq 1$ in $K_n$. Thus, we can assume that $k > 0$ in this case. Now assume that $a_n$ occurs in $\tau _j$ for some $1\leq j\leq r$. The group $K_n$ is isomorphic to $F\ast K$ where $F$ is the free subgroup generated by $\{a_i : i\leq n-3\}$ and $K$ is the subgroup generated by $\{a_{n-2}, a_{n-1}, a_{n}\}$; notice that $K\cong \langle x, y, z \ | x^2y^{-1}[x,y]^qz^2 = 1\rangle$ with isomorphism $a_{n-2}\to x, a_{n-1}\to y, a_{n}\to z$. On the other hand, the group $K_n$ has a quotient isomorphic to $F\ast G\ast \Z _2$ with a quotient map $\phi :K_n\to F\ast G\ast \Z _2$ where $F$ is a free group on the set $\{\phi (a_i) : i < n-2\}, G\cong \langle x, y \ | \ x^{2}y^{-1}[x,y]^q \rangle $ and $G$ is generated by $\{\phi (a_{n-2}), \phi (a_{n-1})\}$ with isomorphism $\phi (a_{n-2})\to x, \phi (a_{n-1})\to y$, and $\phi (a_n)$ generates $\Z _2$. (We abuse the notation by denoting the image of $F$ again by $F$ since $\phi |_{F}$ is an isomorphism.) For all $l\geq 0$, let $J_l = \{j : a_{n-l} \ \mathrm{occurs \ in} \ \tau _j, \ \mathrm{but} \ a_i \ \mathrm{does \ not \ for \ all} \ i > n-l\}$. Notice that the set $J_l, l\geq 0$ are pairwise disjoint.

  \medskip

   On the other hand, by induction on $r$ in $\mathcal{F}$, we may assume that if $W = 1$ in $K_n$, then none of $\tau _i, V_i, 1\leq i\leq r$ contains a letter $a_i$ for $i\leq n-3$. Indeed, let us assume the opposite. Then, considering the quotient $K_n\to F\oplus K$, we reduce it to the case that for some $j\leq n-3$ and $i\in \{1, \dots , r\}$, $a_j$ occurs in $\tau _{i}$. Let $s = \min \{j : a_j \ \mathrm{occurs \ in} \  \tau _i \ \mathrm{for \ some} \ 1\leq i\leq r\} $ and let $a_s$ occur in $\tau _l$ for some $1\leq l\leq r$. We may also assume that $a_j$ does not occur in $\tau _i, V_i, 1\leq i\leq r$ for any $j < s$. Let now $\sigma $ be the sum of exponents of $a_s$ in $\tau _l$. Then for all $1\leq i\leq r$, the sum of exponents of $a_{s(i)}$ equals $\sigma $ where $s(i) = \min \{j : a_j \ \mathrm{occurs \ in} \  \tau _i\}$. Without loss of generality, we may assume that $\sigma \geq 0$. If $\sigma > 0$ then it is straightforward to embed $K_n$ in $\mathrm{Homeo}_+(\R)$ such that $\tau _i(x)\geq x$ for all $1\leq i\leq r, x\in \R$ and $\tau _l(x) > x$ for all $x\in \R$. This implies that $W\neq 1$ in $K_n$ and we are done. If $\sigma = 0$, then for all $s < i \leq n$ we impose relations $a_sa_ia_s^{-1} = a_i^{p_{i}}$ in $K_n$ where $p_{i}$ are positive integers; let $K_n'$ be this quotient. By choosing $1<<p_{s+1} << \dots << p_{n}$,  we see that the subgroup of $K_n'$ generated by $a_{s+1}, \dots, a_n$ is word-hyperbolic and considering its action at the boundary we obtain $W\neq 1$ in $K_n'$ hence $W\neq 1$ in $K_n$.

\medskip 

   Thus, all $\tau _i, V_i, 1\leq i\leq r$ are words in the alphabet $\{a_{n-2}, a_{n-1}, a_n\}$. Then $J_l = \emptyset $ for all $l\geq 3$, and since $\tau _i, 1\leq i\leq r$ is not a power of any letter $a_j$ \footnote{If $\tau _i = a_k^r$ for some $r\in \Z\backslash \{0\}$, then its image in the Abelianization under the map $\mathcal{F}\to \mathcal{F}/[\mathcal{F}, \mathcal{F}]\cong \displaystyle \mathop{\oplus}_{i\in \Z}\Z$ will be represented with a vector $re_j$ for some $j\in \{n-2, n-1, n\}$ where $e_k = \langle \dots , 0, 0, 0, 1, 0, 0, \dots \rangle $ whith ``1'' stands at the $k$-th coordinate. Then a non-zero linear combination of $re_{n-2}, re_{n-1}, re_{n}$ with non-negative integral coefficients must be collinear with a vector $\langle \dots , 0, 0, 2, -1, 2, 0, \dots \rangle$ but this is impossible.}, we can claim that $J_l = \emptyset $ for all $l\geq 2$, moreover, for all $i\in J_0$, $\tau _i$, as a reduced word, contains all three letters $a_n, a_{n-1}, a_{n-2}$ or is a reduced word in the alphabet $\{a_n, a_{n-1}\}$ that contains both letters. In the former case, we have $J_l = \emptyset $ for all $l\neq 0$ hence $K_n = K_n/\langle \tau _i, i\in \cup _{l\geq 1}J_l\rangle $. Then, $\tau _i$ is a generalized torsion in the group $K \cong \langle x, y, z \ | x^2y^{-1}[x,y]^qz^2 = 1\rangle$,  However, the latter is a bi-orderable group by a result of John Labute \cite{JL}. \footnote{It is shown that if $G = \langle x_1, \dots , x_m, y_1, \dots , y_n \ | \ u = v\rangle $ where $u\in A  = F(x_1, \dots , x_m), v\neq 1$ and for some $d$, $u\in \gamma _d(A)$, is not a proper power modulo $\gamma _{d+1}(A)$ then $G$ is residually torsion-free nilpotent, hence bi-orderable. We can take $A = F(x,y), B = F(z), u=x^2y^{-1}[x,y]^q\in A, v=z^{-2}\in B$. Then $\gamma _2(A) = [A, A], u\equiv x^2y^-1 (\mathrm{mod}(\gamma _(A))$. Since $A/\gamma_2(A) \cong \mathbb{Z}^2$, the image of $u$ is the vector $(2,-1)$ which is primitive. Hence $x^2y^{-1}$ is not a proper power in $A/\gamma_2(A)$, therefore $u$ is not a proper power modulo $\gamma _2(A)$. Also, $v = z^{-2}\neq 1$}   So, we will assume that for all $i\in J_0$, $\tau _i$ is a reduced word in the alphabet $\{a_n, a_{n-1}\}$ that contains both letters.  Then for all $i\in J_1$, $\tau _i$ is a reduced word in the alphabet $\{a_{n-1}, a_{n-2}\}$ that contains both the letters $a_{n-1}$ and $a_{n-2}$. We choose any $i\in J_0$ and let $\tau =\tau (u,v)$ be the reduced word in the free group $\langle u, v\rangle \cong \F _2 $ with a chosen standard generating set such that the word $\tau _i = \tau _i(a_{n-1},a_n)$ is represented with the word $\tau $.   

   \medskip 

    Using the quotient map $\mathcal{F}\to \mathcal{F}/[\mathcal{F}, \mathcal{F}]\to \R$, we also see that for all $i\in J_0$, the sum of exponents of $a_n$ and of $a_{n-1}$ are both zero. Then for all $i\in J_1$, the sum of exponents of $a_{n-1}$ and $a_{n-2}$ is zero. Thus, the word $\tau $ belongs to the commutator subgroup $\langle u, v\rangle ^{(1)}$ of the free group $\langle u, v\rangle$. Then the image of $\tau(y,z)$ will be the identity element in the quotient $K' = \langle x, y, z \ | x^2y^{-1}[x,y]^qz^2 = 1\rangle/\langle [y,z]=1\rangle $. On the other hand, $K'$ has a quotient $K''(m) = \langle x, y, z \ | x^2y^{-1}[x,y]^qz^2 = 1\rangle/\langle z=y^m\rangle \cong \langle x, y \ | x^2y^{-1}[x,y]^qy^{2m} = 1\rangle $ for all $m\geq 1$. For any given $R > 0$, if $m$ is sufficiently large, $K''(m)$ is a word hyperbolic group and there will be no relation of length less than $R$ between $x$ and $y$ in the group $K''(m)$; hence the image of $\tau(x,y)$ and $W$ in $K''(m)$ is not the identity element. Then the group $K_n/\langle \tau _i, i\in J_0\cup \displaystyle \mathop{\cup }_{l\geq 2}J_l\rangle $ has a quotient $K_n/\langle \tau _i, i\in  \displaystyle \mathop{\cup }_{l\geq 2}J_l, a_n=a_{n-1}^m\rangle $ and for sufficiently large $m$ and for all $i\in J_1$, the image of $\tau _i$ and $W$ are not the identity elements. Contradiction.

\vspace{1cm} 

 {\em Acknowledgement:} We are very grateful to an anonymous referee for pointing out two counterexamples to a previous version of this paper, which helped us quickly identify a flaw (omissions) in our argument. These counterexamples were instrumental to our understanding.

 \end{document}